\documentclass[11pt]{amsart}

\usepackage{amssymb,amscd,amsthm,amsxtra}
\usepackage{latexsym}
\newtheorem{proposition}{Proposition}
\newtheorem{theorem}{Theorem}
\newtheorem{lemma}{Lemma}
\begin{document}

\title
[Mass concentration for the quintic NLS in 1d]
{\bf Mass concentration phenomenon for the quintic nonlinear
Schr\"odinger equation in 1d}

\author{Nikolaos Tzirakis}

\address{Nikolaos Tzirakis\\
Department of Mathematics\\
University of Toronto\\
Toronto, Ontario, Canada, M5S 2E4\\}

\email{tzirakis@math.toronto.edu}

\date{\today}



\thanks{This material is based upon work supported by the National Science Foundation under agreement No. DMS-0111298. Any opinions, findings and conclusions
 or recommendations expressed in this material are those of the author and do not reflect the views of the National Science Foundation.}
\begin{abstract}
We consider the $L^{2}$-critical quintic focusing nonlinear Schr\"odinger equation (NLS) on ${\bf R}$. It is well known that 
$H^{1}$ solutions of the aforementioned equation blow-up in finite time. In higher dimensions, for $H^{1}$ spherically symmetric blow-up solutions
 of the $L^{2}$-critical focusing NLS, there is a minimal amount of concentration of the $L^{2}$-norm (the mass of the ground state) at the origin.
 In this paper we prove the existence of a similar phenomenon for the 1d case and rougher initial data, $(u_{0}\in H^{s},\ s<1)$, without any 
additional assumption.
\end{abstract}
\maketitle 
\section{Introduction}
This paper continues the investigation of the quintic nonlinear Schr\"odinger equation (NLS) in one dimension that we started in \cite{tz}. 
\begin{align}\begin{split}
&iu_{t}+u_{xx}\pm |u|^{4}u=0\\
&u(x,0)=u_{0}(x)\in H^{s}({\bf R}),t\in {\bf R}.
\end{split}
\end{align}
The $(+)$ sign in front of the nonlinearity corresponds to the focusing NLS while the $(-)$ sign to the defocusing. The Cauchy problem
for equation $(1)$ is known to be locally well-posed in $H^{s}({\bf R})$ for $s>0$, \cite{cw}. A local result also exists for $s=0$, but
 the time of existence depends on the profile of the data as well as the norm. NLS is an infinite dimensional Hamiltonian system with 
energy space $H^{1}$. It also has a scaling property. Thus
 $u(x,t)$ is a solution of $(1)$ with initial data $u_{0}$ if and only if 
$$u^{\lambda}(x,t)=\frac{1}{\lambda^{1/2}}u(\frac{x}{\lambda},\frac{t}{\lambda^{2}})$$ is a solution to the same equation 
with initial data $u_{0}(\frac{x}{\lambda})$. In \cite{tz} we extend the local existence theorem for the defocusing NLS for all times. We did
so by iterating the local result in the appropriate norms. To iterate the local 
result by standard limiting arguments we just need an apriori bound for our solutions in $H^{s}$. 
This bound comes from the next theorem that we proved in \cite{tz}:
\begin{theorem}
Let $u$ be a global $H^{1}$ solution to $(1)$ with the $(-)$ sign. Then for any $T>0$ and $s>4/9$ we have that 
$$\sup_{0 \leq t \leq T}\|u(t)\|_{H^{s}} \lesssim C_{(\|u_{0}\|_{H^{s},T})}$$
\\
where the right hand side does not depend on the $H^{1}$ norm of $u$. 
\end{theorem}
{\bf Remark.} Note that in the focussing case (where in front of the
 nonlinearity we have the plus instead of the minus sign) we can also proved global well-posedness for $4/9<s\leq 1/2$, but
 with the crucial assumption that $\|u_{0}\|_{L^{2}}<\|Q\|_{L^{2}}$ where Q is the unique positive solution (up to translations) of
$$Q_{xx}-Q+|Q|^{4}Q=0.$$
In \cite{mw1}, Q was solved explicity as $Q(x)=3^{\frac{1}{4}}/\sqrt{\cosh (2x)}$ and then $\|Q\|_{L^{2}}^{2}=\frac{\sqrt{3}\pi}{2}$. In
 the same paper it was also proved a result that we will use below, namely that $C=\|Q\|_{L^{2}}^{-4}$ is the best constant in the Gagliardo-Nirenberg inequality
$$\frac{1}{6}\|u\|_{L^{6}}^{6} \leq \frac{C}{2}\|\nabla u\|_{L^{2}}^{2}\|u\|_{L^{2}}^{4}.$$
We used the ``I-method'' that was recently introduced by J. Colliander, M. Keel, G. Stafillani, H.Takaoka, and T.Tao, \cite{ckstt1,ckstt3,ckstt4,ckstt5}. 
This method allows us to define a modification of the energy functional, that is ``almost conserved'' that is, its time derivative decays with respect to
 a very large parameter. Since an implementation of this method gives the main result of this paper also, the details of the method are delayed until 
the next section. As we mentioned above for the focusing case, the solutions blows up in $H^{1}$, in finite time. An elementary proof of the existence of blow-up
 solutions is known since the 60's, but is based on energy constraints and is not constructive, \cite{ss}. In particular no qualitative information
 of any type of the blow-up dynamics is obtained. A lower estimate for the blow-up solutions in $H^{1}$ is given by Theorem 2 below using the 
scaling and the local existence theorem, \cite{tc}.
\begin{theorem}
Let $[0,T^{\star})$ is the maximal interval of existence of the following Cauchy problem:
\begin{align}\begin{split}
&iu_{t}+u_{xx}+|u|^{4}u=0\\
&u(x,0)=u_{0}(x)\in H^{1}({\bf R}),t\in {\bf R}.
\end{split}
\end{align}
If $u_{0} \in H^{1}$ is such that $T^{\star}<\infty$, then there exists a $C$ such that
$$\|u_{x}\|_{L^{2}}\geq \frac{C}{(T^{\star}-t)^{\frac{1}{2}}}$$
for $0\leq t<T^{\star}$.
\end{theorem}
This bound is often called the scaling bound. It is also fairly easy to show that $\|u(t)\|_{L^{p}}$ blows up for $p>2$. In particular we have
\\
$$\|u\|_{L^{p}}\geq \frac{C}{(T^{\star}-t)^{\frac{1}{4}-\frac{1}{2p}}}$$
for $p>2$.
\\
\\
Because it is related to the scaling symmetry of the problem, the above lower bound has long been conjectured to be optimal. But in 1988 Landman, 
Papanicolaou, Sulem, Sulem, \cite{lpss}, suggested that the correct and stable blow up speed is a slight correction to the scaling bound:
\\
$$\|u_{x}\|_{L^{2}}\sim \sqrt {\frac{\log |\log |T^{\star}-t||}{T^{\star}-t}}$$
\\
In this frame Perelman, \cite{pg}, has constructed a family of blowing up solutions for which
\\
$$\left( \frac{\log |\log |T^{\star}-t||}{T^{\star}-t}\right)^{-\frac{1}{4}}\|u(t)\|_{L^{\infty}}\rightarrow c>0$$
\\
as $t\rightarrow T^{\star}$, which is very close but different than the scaling bound. Moreover for initial data in some special class, Merle and Raphael, 
\cite{mr,mr1}, recently showed that for $t$ close to $T^{\star}$, there is a universal constant $C^{\star}$ such that
\\
$$\|u_{x}\|_{L^{2}}\leq C^{\star}\sqrt {\frac{\log |\log |T^{\star}-t||}{T^{\star}-t}}$$
\\
as suggested by the numerics in \cite{lpss}. Finally it is worth noting that an easy application of the pseudoconformal transformation
 yields interesting information on the blow up solutions. In particular we can show that some solutions blow-up twice as fast as
 the scaling bound. For details see \cite{bw} and \cite{mw}. The above results show in particular that at least two different blow up estimates are actually achieved.
\\
\\
Another property of the blow up solutions in the critical case is the phenomenon of mass concentration, \cite{tc} and \cite{ss}. For
 $H^{1}$ solutions, there is a concentration of a finite amount of mass in a neighborhood of the focus of width sligtly larger than $(T^{\star}-t)^{1/2}$.
 For radial initial data in dimension $d \geq 2$ there is a precise lower bound on the amount of concentrated mass in terms of the mass of the
 ground state Q, \cite{mt}. More precisely we have:
\\
\\
$\bullet$ Let $d\geq 2$ and $\gamma:(0,\infty) \rightarrow (0,\infty)$ be any function such that $\gamma(s) \rightarrow \infty$  and
 $s^{1/2}\gamma(s) \rightarrow 0$ as $s \downarrow 0$. Finally let $u_{0} \in H^{1}({\bf R^{d}})$ is radial symmetric. Then if $u(x,t)$ is the maximal solution
of the equivalent of $(2)$ in higher dimensions and $T^{\star}<\infty$ we have
$$\liminf_{t\uparrow T^{\star}} \|u(t)\|_{L_{\{|x|<|T^{\star}-t|^{1/2}\gamma(T^{\star}-t)\}}^{2}} \geq \|Q\|_{L^{2}}.$$
\\
where $Q$ is the ground state solution of the elliptic equation $Q_{xx}-Q+|Q|^{4}Q=0$.
\\
\\
In the nonradial case and in dimension $d=1$ this was generalized by Nawa, \cite{hn}, using concentration compactness techniques, \cite{l1,l2}. In addition
to the scaling properties of the NLS equation, the main ingredients in the proof that $H^{1}$ blowup solutions concentrate at least the mass
 of the ground state are:
\\
\\ 
i) The conservation of mass
$$\|u(t)\|_{L^{2}}=\|u_{0}\|_{L^{2}}$$
and the energy
$$E(u)(t)=E(u_{0})$$ 
where $$E(u)=\frac{1}{2}\int{|u_{x}(t)|^{2}dx}-\frac{1}{6} \int{|u(t)|^{6}dx}$$
\\
ii) a precise Galiardo-Nirenberg inequality which implies that nonzero $H^{1}-$functions of non-positive energy have at least ground state mass. 
\\
\\
The purpose of this paper is to investigate the mass concentration phenomenon in $H^{s}$, for $s<1$, where the conservation of energy
 cannot be used. Using the $I-method$ we show that solutions of $(2)$ with a finite maximal (forward) existence interval are
 expected to concentrate at least the $L^{2}$ mass of the ground state in $H^{s}$, for $s<1$. More precisely we have:
\begin{theorem}
Suppose  $H^{s} \ni u_{0}\longmapsto u(t)$ with $s>0$ solves $(2)$ on
 the maximal interval of existence $[0,T^{\star})$ with $T^{\star}<\infty$. Then for any $1>s>\frac{10}{11}$ 
there exists a positive function $\gamma(x)\uparrow \infty$ arbitrarily slowly
 as $x\downarrow 0$ and a real function $z(t)$ such that
$$\limsup_{t\uparrow T^{\star}} \|u(t)\|_{L_{\{|x-z(t)|<(T^{\star}-t)^{\frac{s}{2}}\gamma(T^{\star}-t)\}}^{2}} \geq \|Q\|_{L^{2}}$$
\end{theorem}
${\bf Remark\ 1.}$ In a recent preprint, J. Colliander, S. Raynor, C. Sulem and J .D. Wright, \cite{crsw} consider the $2d$ focusing critical NLS
and proved a similar theorem with the additional assumption of radial symmetry. The radial symmetry assumption is needed in order to pass from weak
 to strong convergence since the general embedding $H^{1}({\bf R^{d}})\hookrightarrow L^{2}({\bf R^{d}})$ is not compact. But as the four authors
 note in \cite{crsw}, one can utilize the concentration compactness method of Lions, \cite{l1,l2} and prove the analogus theorem in 2d. 
The $1d$ case that we are dealing with, have some similar
 features but also significant differences. First of all in the 1d case the radial assumption doesn't play a role. More precisely in 1d, radial symmetry
 is not enough for a bounded sequence in $H^{1}$ to have a strongly convergent subsequence in $L^{p}$ for $2<p<\infty$ although the latter is true if
 we further assume that the sequence in question is a nonincreasing function of $|x|$ for every $n \geq 0$.(For the above discussion the reader can also consult
 \cite{ws}). So we have to prove Theorem 3 by implementing different techniches than the techniques used in \cite{crsw}. 
Second the nonlinearity has a fifth power,
 and thus the correction terms in the ``modified energy''  have larger growth. We take advantage of the fact that at each step we work on $[0,\delta]$ and we
prove a stronger proposition about the decay of the ``modified energy'' and thus somehow we balance the additional correction terms with the greater decay that we prove.
 Finally the crucial Lemma 3 that we use in 1d, is true only if the frequencies of the two solutions are seperated. In higher dimensions
 the analogous Lemma holds in general, \cite{bo}, although we avoid this difficulty in 1d by analyzing further the correction terms of the ``modified energy'',
 see Proposition $5$.
\\
\\
${\bf Remark\ 2.}$ As we mentioned before, to prove the theorem we use a combination of the concentration compactness and the $I-method$. Since
the energy is infinite for initial data in $H^{s}$ we define a ``modified energy'', $E(Iu)$ which is finite, where $I:H^{s}\rightarrow H^{1}$
 is a multiplier operator defined below. The crucial step is to prove that the modified total energy grows more slowly than the modified kinetic energy
 $$\frac{1}{2}\int |Iu_{x}(t)|^{2}dx.$$
These two steps are shown in 
Propositions $5$ and $3$ respectively. Note that Proposition 5 relys on the local theory that we shall establish in Proposition 1. 
\\
\\
${\bf Remark\ 3.}$ Let $p(s)$ be a number that depends on $s$ and for the range of $s$ in Theorem 3, ($10/11<s<1$), it is $p(s)<2$. 
The statement that the modified total energy grows more slowly than the modified kinetic energy is reflected exactly on $p(s)<2$
 and is proven in Proposition 3 below. Note
also that our concentration width $(T^{\star}-t)^{\frac{s}{2}}$ is larger than $(T^{\star}-t)^{\frac{1}{2}}$ with which ground state mass
 concentration is conjectured to occur.
\\
\\
Two quick by-products of the above theorem are the following. The first is the conjecture that tiny $L^{2}$ mass concentration cannot occur when
 $u_{0}\in L^{2}$, a question that was asked in \cite{mv}. See also the relevant result of J. Bourgain, \cite{bo}. The second is the following
 lemma which as we mention on the first page is basically a result of the work in \cite{tz}. 
\begin{lemma}
If $u_{0}\in H^{s}$, $s>\frac{10}{11}$ and $\|u_{0}\|_{L^{2}}<\|Q\|_{L^{2}}$, then the initial value problem $(2)$ is globally well-posed.
\end{lemma}
We end in this section by introducong some useful notation. In what follows we use $A \lesssim B$ to denote an estimate of the form $A\leq CB$ for some constant $C$.
If there exist constants $C$ and $D$ such that $DB \leq A \leq CB$ we
say that $A \sim B$, and $A \ll B$ to denote an estimate of the form $A \leq cB$ for small constant $c>0$. In addition $\langle a \rangle:=1+|a|$ and
$a\pm:=a\pm \epsilon$.

\section{Linear and Bilinear Estimates}

Before we state the linear and bilinear estimates that we will use throuout this 
paper we recall some basic facts about the $X^{s,b}$ spaces. 
For an equation of the form
\begin{equation}
iu_{t}-\phi(-i \nabla)u=0
\end{equation}
where $\phi$ is a measurable function let $X^{s,b}$ be the completion of $\Bbb{S}(\Bbb R^{d+1})$
with respect to 
$$\|u\|_{X^{s,b}}=\|\langle \xi \rangle^{s}\langle \tau+\phi(\xi)\rangle^{b}\hat{u}(\xi,\tau)\|_{L_{\xi}^{2}L_{\tau}^{2}}.$$
From the above definition it is clear that the dual space of $X_{\tau==\phi(\xi)}^{s,b}$
is $X_{-\tau=-\phi(\xi)}^{-s,-b}$. Furthermore for a given interval $I$, we define 
$$\|f\|_{X^{s,b}(I)}=\inf_{\tilde{f}_{\vert I}=f} \|\tilde{f}\|_{X^{s,b}}$$
In our case, the interval of existence of the local solutions will be $[0,\delta]$ and we write $X_{\delta}^{s,b}=X_{[0,\delta]}^{s,b}$. 
Since conjugate solutions won't play any role in our arguments
from now on we ommit any reference to the difference between $u$ and $\bar{u}$.
We know that if $u$ is a solution of $(3)$ with $u(0)=f$ and $\psi$ is a cut-off function in $C_{0}^{\infty}$
with support of ${\psi} \subset (-2,2)$, $\psi=1$ on $[0,1]$, $\psi(-t)=\psi(t)$, 
$\psi(t) \geq 0$, $\psi_{\delta}(t)=\psi(\frac{t}{\delta})$ then if $0<\delta\leq 1$ we have that for $b\geq 0$:
\begin{equation}
\|\psi_{1}u\|_{X^{s,b}} \leq C\|f\|_{H^{s}}
\end{equation}
In addition if $\nu$ is a solution of
$$i\nu_{t}-\phi(-i \nabla)\nu=F$$ 
\\
with $\nu(0)=0$ then for $b^{'}+1 \geq b \geq 0 \geq b^{'}>-\frac{1}{2}$:
\begin{equation}
\|\psi_{\delta} \nu\|_{X^{s,b}} \leq C\delta^{1+b^{'}-b}\|F\|_{X^{s,b^{'}}}.
\end{equation}
The proofs of $(4)$ and $(5)$ can be found in \cite{gtv}. 
The Strichartz estimates for the Schr\"odinger equation on $\Bbb R^{d}$ state that 
for $q,r \geq 2$ such that $(d,q)\ne (2,2)$ and $0 \leq \frac{2}{q}=d(\frac{1}{2}- \frac{1}{r})<1$, we have that
\begin{equation}
\|e^{it\Delta}u_{0}\|_{L_{t}^{q}L_{x}^{r}} \lesssim \|u_{0}\|_{L^{2}(\Bbb R^{d})}.
\end{equation}
In particular, in $1d$ we have:
$$\|e^{it\partial_{x}^{2}}u_{0}\|_{L_{t}^{6}L_{x}^{6}} \lesssim \|u_{0}\|_{L^{2}(\Bbb R)}$$
and
$$\|e^{it\partial_{x}^{2}}u_{0}\|_{L_{t}^{\infty}L_{x}^{2}} \lesssim \|u_{0}\|_{L^{2}(\Bbb R)}$$
which by a standard argument gives 
\begin{equation}
\|u\|_{L_{t}^{6}L_{x}^{6}} \lesssim \|u\|_{X_{\delta}^{0,1/2+}}
\end{equation}
and 
\begin{equation}
\|u\|_{L_{t}^{\infty}L_{x}^{2}} \lesssim \|u\|_{X_{\delta}^{0,1/2+}}.
\end{equation}
By Sobolev embedding theorem in $1d$, $(8)$ implies that
\begin{equation}
\|u\|_{L_{t}^{ \infty}L_{x}^{ \infty}} \lesssim \|u\|_{X_{\delta}^{1/2+,1/2+}}.
\end{equation}
Also by interpolation between $(7)$ and the trivial estimate,
\begin{equation}
\|u\|_{L_{t}^{2}L_{x}^{2}}=\|u\|_{X_{\delta}^{0,0}}
\end{equation}
we get 
\begin{equation}
\|u\|_{L_{t}^{p}L_{x}^{p}} \lesssim \|u\|_{X_{\delta}^{0,(1/2+)\cdot(\frac{3}{2}-\frac{3}{p})}}
\end{equation}
for any $2 \leq p \leq 6$.
\\
\\
The dual version of $(6)$ gives
\begin{equation}
\|u\|_{X_{\delta}^{0,-1/2-}} \lesssim \|u\|_{L_{t}^{q^{'}}L_{x}^{r^{'}}}
\end{equation}
where $r^{'}$ and $q^{'}$ are the conjugate exponents of $r$ and $q$ respectively.
Interpolation with the trivial estimate
\begin{equation}
\|u\|_{X_{\delta}^{0,0}}=\|u\|_{L_{t}^{2}L_{x}^{2}}
\end{equation}
gives that
\begin{equation}
\|u\|_{X_{\delta}^{0,-1/2+}} \lesssim \|u\|_{L_{t}^{q^{'}+}L_{x}^{r^{'}+}}
\end{equation}
and also that
\begin{equation}
\|u\|_{X_{\delta}^{1,-1/2+}} \lesssim \|u\|_{L_{t}^{q^{'}+}W_{x}^{1,r^{'}+}}
\end{equation}
for any $\frac{2}{q}+\frac{1}{r}=\frac{1}{2}$.
\\
\\
As we state in the introduction we prove that the modified energy grows more slowly than the modified kinetic energy in Proposition 3. We can take advantage of the fact
 that we work on a small interval $[0,\delta]$ and improve the decay of the ``modified energy''. To do so we state the following Lemma that we can
 find in \cite{gtv} and \cite{hp}.
\begin{lemma}
If $1/2>b>b^{\prime}\geq 0$ and $s \in {\bf R}$ then the following embedding is true:
$$\|f\|_{X_{\delta}^{s,b^{\prime}}} \lesssim \delta^{b-b^{\prime}}\|f\|_{X_{\delta}^{s,b}}$$
\end{lemma}
The second lemma that we state in this section is an improved bilinear Strichartz type estimate. It is due to Bourgain, \cite{bo}. As we mentioned before
 a general analog of Lemma 3 holds for $d \geq 2$, see for example \cite{ckstt6}.
\begin{lemma}
Let $u$ and $v$ be any two Schwartz functions whose support of Fourier transform is in $|\xi| \sim M$ and
$| \xi| \ll M$ respectively and $M \gg 1$.Then
$$\|(D_{x}^{ \frac{1}{2}}u)v\|_{L_{t}^{2}L_{x}^{2}}=\|(D_{x}^{ \frac{1}{2}} \bar{u})v\|_{L_{t}^{2}L_{x}^{2}}
 \lesssim \|u\|_{X^{0,1/2+}} \|v\|_{X^{0,1/2+}}.$$
\end{lemma}
\section{The I-method and the proof of Theorem 3}
As we mentioned above, the basic step towards Theorem 3 is the fact that the ``modified total energy'' decays more slowly than the ``modified kinetic energy''. 
To prove the last statement we iterate the local solutions for the new modified system
\begin{align}\begin{split}
&iIu_{t}+Iu_{xx}+I(|u|^{4}u)=0\\
&Iu(x,0)=Iu_{0}(x)\in H^{1}({\bf R}),t\in {\bf R}.
\end{split}
\end{align}
So let's define the I-operator. We introduce as in  \cite{ckstt1,ckstt4}, a radial $C^{\infty}$, monotone multiplier, taking values in [0,1], 
 where:
\[m(\xi):= \left\{\begin{array}{ll}
1 & \mbox{if $|\xi|<N$}\\
(\frac{|\xi|}{N})^{s-1} & \mbox{if $|\xi|>2N$}
\end{array}
\right.\]
and we define $I:H^{s} \rightarrow H^{1}$ by $\widehat{Iu}(\xi)=m(\xi)\hat{u}(\xi).$ The operator $I$ is smoothing
of order $1-s$ and we have that:
\begin{equation}
\|u\|_{X_{\delta}^{s_{0},b_{0}}} \lesssim \|Iu\|_{X_{\delta}^{s_{0}+1-s,b_{0}}} \lesssim N^{1-s}\|u\|_{X_{\delta}^{s_{0},b_{0}}}
\end{equation}
for any $s_{0},b_{0}\in {\bf R}$.
\\
\\
$\bf{ Remark.}$ It is shown in \cite{ckstt4} that if
$$\|uv\|_{X^{s,b-1}} \lesssim \|u\|_{X^{s,b}} \|v\|_{X^{s,b}}$$
then
$$\|I(uv)\|_{X^{1,b-1}} \lesssim  \|Iu\|_{X^{1,b}} \|Iv\|_{X^{1,b}}$$
\\
where the constants in the above inequality are intependent of
$N$. From now on we use this fact and refer to it as the
``interpolation lemma''. For details see \cite{ckstt4}.
\begin{proposition}
Let $s>10/11$ and consider the equation
\begin{align}\begin{split}
&iIu_{t}+ (Iu)_{xx}+I(|u|^{4}u)=0
\end{split}
\end{align}
\\
with initial data $Iu(x,0)=Iu_{0}$. Then there exists a 
$$\delta \sim (\|Iu_{0}\|_{H^{1}})^{-4-\epsilon}$$
such that for all times in $[0,\delta]$, the above problem is locally well-posed and
$$\|Iu\|_{X_{\delta}^{1,1/2+}} \lesssim \|Iu_{0}\|_{H^{1}}.$$
\end{proposition}
\begin{proof}
By Duhamel's formula the equation $(18)$ is equivalent to
$$Iu(t)=\psi_{1}(t)e^{it\partial_{x}^{2}}(Iu_{0})+i\psi_{ \delta}(t) 
\int_{0}^{t}e^{i(t-s)\partial_{x}^{2}}I(|u|^{4}u)(s)ds$$
By $(4)$ and $(5)$ and the fact that $\delta \leq 1$ we have
$$\|Iu\|_{X_{\delta}^{1,1/2+}} \lesssim \|Iu_{0}\|_{H^{1}}+\|I(|u|^{4}u)\|_{X_{\delta}^{1,-1/2+}}$$
Now recall the dual Strichartz estimate, equation $(15)$
$$\|u\|_{X_{\delta}^{1,-1/2+}} \lesssim \|u\|_{L_{t}^{q^{\prime}+}W_{x}^{1,r^{\prime}+}}$$
where $\frac{2}{q}=\frac{1}{2}-\frac{1}{r}$. Thus for $r^{\prime}=2-$ we have
$$\|I(|u|^{4}u)\|_{X_{\delta}^{1,-1/2+}} \lesssim \|I(|u|^{4}u)\|_{L_{t}^{1+\epsilon}H_{x}^{1}} \lesssim \delta^{1-\epsilon}\|I(|u|^{4}u)\|_{L_{t}^{\infty}H_{x}^{1}}.$$
Since for $s>1/2$, $H^{s}$ is a Banach algebra we have that
$$\||u|^{4}u\|_{H_{x}^{s}} \lesssim \|u\|_{H_{x}^{s}}^{5}$$
which by the interpolation lemma quickly translates to
$$\|I(|u|^{4}u)\|_{H_{x}^{1}} \lesssim \|Iu\|_{H_{x}^{1}}^{5}.$$
But then
$$\|Iu\|_{X_{\delta}^{1,1/2+}} \lesssim \|Iu_{0}\|_{H^{1}}+\delta^{1-\epsilon}\|Iu\|_{L_{t}^{\infty}H_{x}^{1}}^{5} 
\lesssim \|Iu_{0}\|_{H^{1}}+\delta^{1-\epsilon}\|Iu\|_{X^{1,1/2+}}^{5}$$
and by standard iteration arguments, see \cite{kpv}, we have that the system is locally well-posed for
$$\delta^{1-\epsilon}\|Iu_{0}\|_{H^{1}}^{4}<\frac{1}{2}.$$
\end{proof}
We also need an analog of Theorem 2 for the I-system $(16)$. Let $\nabla ^{s}$ denote the operator which on the Fourier side is given by 
$\widehat{\nabla^{s}u}(\xi)=|\xi|^{s}\hat{u}(\xi)$. It then follows by the definition of the Japanese bracket that on the Fourier side the 
$<\nabla>$u is given by $(1+|\xi|)\hat{u}(\xi)$.
\begin{proposition}
If $H^{s} \ni u_{0}\longmapsto u(t)$ with $s>10/11$ solves $(2)$ for all t close enough to $T^{\star}$ in the maximal finite interval
 of existence $[0,T^{\star})$ then
$$\|I\langle \nabla \rangle u(t)\|_{L^{2}} \geq C(T^{\star}-t)^{-\frac{s}{2}}$$
\end{proposition}
\begin{proof}
Since we know that 
$$\|I\langle \nabla \rangle u(t)\|_{L^{2}} \geq \|u(t)\|_{H^{s}}$$ 
it suffices to show that 
$$\|\nabla^{s}u(t)\|_{L^{2}}\geq C(T^{\star}-t)^{-\frac{s}{2}}.$$
\\
We assume that $\|\nabla^{s}u(t)\|_{L^{2}}>1$ since otherwise we can change variables to put the time origin near to $T^{\star}$. Now fix $t \in [0,T^{\star})$
 and consider
$$v^{t}(\tau,x)=\lambda^{-1/2}u(t+\frac{\tau}{\lambda^{2}},\frac{x}{\lambda})$$
where $\lambda=\|\nabla^{s}u(t)\|_{L^{2}}^{\frac{1}{s}}$. By scaling invariance $v^{t}(\tau,x)$ is a solution to $(2)$. Moreover an easy calculation shows that
$$\|v^{t}(0,x)\|_{L^{2}}=\|u_{0}\|_{L^{2}}$$
and that
$$\|\nabla^{s}v^{t}(0,x)\|_{L^{2}}=\lambda^{-s}\|u(t,x)\|_{\dot{H}^{s}}=1$$
\\
Thus $\|v^{t}(t,x)\|_{H^{s}}<C$ and by 
the local theory that means that there exists a $\tau_{0}>0$, independent of $t$, such that $v^{t}(t,x)$ is defined on $[0,\tau_{0}]$ and therefore
\\
$$t+\frac{\tau_{0}}{\|\nabla^{s}u(t)\|_{L^{2}}^{\frac{2}{s}}} \leq T^{\star} \Longrightarrow \|\nabla^{s}u(t)\|_{L^{2}}\geq C(T^{\star}-t)^{-\frac{s}{2}}.$$ 
\end{proof}
The last step for the proof of Theorem 3 is the following proposition that for the moment we assume and prove later.
\begin{proposition}
For $s>\frac{10}{11}$ there exists $p(s)<2$ such that the following hold true:
\\
If $H^{s} \ni u_{0}\longmapsto u(t)$ solves $(2)$ on $[0,T^{\star})$ then for all $T<T^{\star}$ there exists $N=N(T)$ such that
$$|E[I_{N(T)}u(T)]| \leq C_{0} \Lambda(T)^{p(s)}$$
with $C_{0}=C_{0}(s,T^{\star},\|u_{0}\|_{H^{s}})$, and $\Lambda(T)$ is given in terms of $N(T)$ by $N(T)=C(\Lambda(T))^{\frac{p(s)}{2(1-s)}}$.
\end{proposition}
 We prove Theorem 3 by using the concentration
 compactness method that was developed by Lions in \cite{l1,l2}.We will need a series of Lemmas. The proof of the first two Lemmas are easy and
 can be found in \cite{tc} on pages 21 and 24 respectively.
\\
\begin{lemma}
Let $u\in L^{2}$ and let the concentration function be defined by
$$\rho(u,t)=\sup_{y\in {\bf R}}\int_{\{|x-y|<t\}}|u(x)|^{2}dx$$
for $t>0$. Then $\rho$ is a nondecreasing function of $t$ and there exists $y(u,t)\in {\bf R}$ such that
\\
$$\rho(u,t)=\int_{\{|x-y(u,t)|<t\}}|u(x)|^{2}dx.$$Moreover if $u \in L^{r}(\Bbb R)$ for some $r>2$, then for all $s,t>0$ and $C=C(r)$ we have
$$|\rho(u,t)-\rho(u,s)| \leq C\|u\|_{L^{r}}^{2}|t-s|^{\frac{r-2}{r}}$$
\end{lemma}
\begin{lemma}
There exists a constant $K$ such that for all $u \in H^{1}$, all $t>0$ and $\rho$ defined above we have
$$\int |u|^{6} \leq K \rho (u,t)^{2}\left( \int |\nabla u|^{2}+t^{-2}\int |u|^{2} \right)$$
\end{lemma}
\begin{lemma}
Let $(u_{n})_{n\geq 0} \subset H^{1}$ be such that
$$\|u_{n}\|_{L^{2}} \leq a <\infty$$
\\
$$\sup_{n\geq 0}\|\nabla u_{n}\|_{L^{2}}<\infty$$
\\
and let $\rho(u_{n},t)$ defined as before. Set 
$$\mu=\lim_{t\rightarrow \infty}\liminf_{n\rightarrow \infty}\rho(u_{n},t).$$
Then there exist a subsequence $(u_{n_{k}})_{n_{k}\geq 0}$,
 a nondecreasing function $\gamma(t)$, and a sequence $t_{k}\rightarrow \infty$ with the following properties:
\\
\\
i) $\rho(u_{n_{k}},.)\rightarrow \gamma(.)\in [0,a]$ as $k\rightarrow \infty$ uniformly on bounded sets of $[0,\infty)$.
\\
\\
ii)$\mu=\lim_{t\rightarrow \infty}\gamma(t)=\lim_{k\rightarrow \infty}\rho(u_{n_{k}},t_{k})=\lim_{k\rightarrow \infty}\rho(u_{n_{k}},t_{k}/2)$. 
\end{lemma}
\begin{proof}
Since $$\mu=\lim_{t\rightarrow \infty}\liminf_{n\rightarrow \infty}\rho(u_{n},t)$$ there exists a $t_{k} \rightarrow \infty$ such that
\begin{equation}
\mu=\lim_{k\rightarrow \infty} \rho(u_{n_{k}},t_{k})
\end{equation}
and thus one part of ii) is evident. To prove the first part note that
$$\rho(u_{n},t) \leq \|u_{n}\|_{L^{2}} \leq a<\infty.$$
In addition since $H^{1}({\Bbb R})\hookrightarrow L^{r}({\Bbb R})$ for some $r$, by the last property of the previous Lemma
$\rho(u_{n},\cdot)$ is H\"older continuous. Therefore i) follows from Ascoli's theorem (after renaming the sequence $n_{k}$). Notice
 that property $(19)$ is still true after passing to a subsequence. For the rest of ii) by $(19)$ and the fact that $\rho(u_{n},\cdot)$ is 
nondecreasing we deduce that
\begin{equation}
\limsup_{k \rightarrow \infty}\rho(u_{n_{k}},\frac{t_{k}}{2}) \leq \limsup_{k \rightarrow \infty}\rho(u_{n_{k}},t_{k})=\mu.
\end{equation}
Next for every $t>0$ we have
$$\liminf_{k \rightarrow \infty} \rho(u_{n_{k}},t) \geq \liminf_{n \rightarrow \infty} \rho(u_{n},t).$$
Now by letting $t \rightarrow \infty$ and using part i) of the Lemma and the definition of $\mu$ we get that
\begin{equation}
\lim_{t \rightarrow \infty} \gamma(t) \geq \mu.
\end{equation} 
Finally, given $t>0$ we have $\frac{t_{k}}{2} >t$ for $k$ large, so that
$$\rho(u_{n_{k}},\frac{t_{k}}{2}) \geq \rho(u_{n_{k}},t)$$
and by letting $k \rightarrow \infty$ by part i) we get
\begin{equation}
\liminf_{k \rightarrow \infty}\rho(u_{n_{k}},\frac{t_{k}}{2}) \geq \mu.
\end{equation}
By (20) and (22) we have that
$$\mu=\lim_{k\rightarrow \infty}\rho(u_{n_{k}},t_{k}/2)$$
Similarly
$$\rho(u_{n_{k}},\frac{t_{k}}{2}) \geq \rho(u_{n_{k}},t)\Rightarrow \sup \rho(u_{n_{k}},\frac{t_{k}}{2}) \geq \rho(u_{n_{k}},t).$$
and by taking $k \rightarrow \infty$ and use (20) we get
\begin{equation}
\mu \geq \lim_{t \rightarrow \infty} \gamma(t).
\end{equation} 
\end{proof}
\begin{lemma}
Let $(u_{n})_{n\geq 0} \subset H^{1}$ be such that
$$\|u_{n}\|_{L^{2}} \leq a <\infty,$$
$$\lim_{n \rightarrow \infty}\|u_{n}\|_{L^{2}}^{2}=b>0$$
and
$$\sup_{n \geq 0}\|\nabla u_{n}\|_{L^{2}}< \infty.$$
Then there exists a subsequence $(u_{n_{k}})_{k \geq 0}$ which satisfies the following:
\\
\\
There exist $(q_{k})_{k\geq 0},(w_{k})_{k\geq 0}\subset H^{1}(\Bbb R)$ such that
\begin{equation}
supp q_{k}\cap supp w_{k}=\emptyset,
\end{equation}
\begin{equation}
|q_{k}|+|w_{k}|\leq |u_{n_{k}}|,
\end{equation}
\begin{equation}
\|q_{k}\|_{H^{1}}+\|w_{k}\|_{H^{1}}\leq C \|u_{n_{k}}\|_{H^{1}},
\end{equation}
\begin{equation}
\lim_{k \rightarrow \infty}\|q_{k}\|_{L^{2}}^{2}=\mu,\ \ \ \ \lim_{k \rightarrow \infty}\|w_{k}\|_{L^{2}}^{2}=b-\mu
\end{equation}
\begin{equation}
\liminf_{k \rightarrow \infty}\{ \int |\nabla u_{n_{k}}|^{2}-\int |\nabla q_{k}|^{2}-\int |\nabla w_{k}|^{2} \} \geq 0,
\end{equation}
\begin{equation}
\lim_{k \rightarrow \infty}\Bigr| \int |u_{n_{k}}|^{p}-\int |q_{k}|^{p}-\int |w_{k}|^{p}\Bigr|=0
\end{equation}
for all $2 \leq p <\infty$.
\end{lemma}
\begin{proof}
We use the sequences $(u_{n_{k}})_{k\geq 0}$ and $(t_{k})_{k \geq 0}$ constructed in the previous Lemma. We fix $\theta ,\phi \in C^{\infty}([0,\infty))$ 
such that $0\leq \theta, \phi \leq 1$ and
\\
\\
$\theta(t)=1$ for $0 \leq t \leq \frac{1}{2}$,\ \ \ \ $\theta(t)=0$ for $t \geq \frac{3}{4}$
\\
\\
$\phi(t)=0$ for $0 \leq t \leq \frac{3}{4}$,\ \ \ \ $\phi(t)=1$ for $t \geq 1$,
\\
\\
and we set 
$$q_{k}=\theta_{k}u_{n_{k}},\ \ \ w_{k}=\phi_{k}u_{n_{k}}$$
where
$$\theta_{k}=\theta \left ( \frac{|x-y(u_{n_{k}},\frac{t_{k}}{2})|}{t_{k}}  \right )\ \ \ \phi_{k}=\phi \left ( \frac{|x-y(u_{n_{k}},\frac{t_{k}}{2})|}{t_{k}}  \right ).$$
Now $(24),(25)$ and $(26)$ are immediate. To prove $(27)$ we estimate
$$\rho(u_{n_{k}},\frac{t_{k}}{2})=\int_{|x-y(u_{n_{k}},\frac{t_{k}}{2})|\leq \frac{t_{k}}{2}} |u_{n_{k}}|^{2} \leq \int |q_{k}|^{2} \leq 
\int_{|x-y(u_{n_{k}},\frac{t_{k}}{2})|\leq t_{k}} |u_{n_{k}}|^{2}$$
$$\leq \int_{|x-y(u_{n_{k}},t_{k})|\leq t_{k}} |u_{n_{k}}|^{2} \leq \rho(u_{n_{k}},t_{k}).$$
Applying the second part of Lemma 6 we immediatelly get
\begin{equation}
\lim_{k \rightarrow \infty}\|q_{k}\|_{L^{2}}^{2}=\mu.
\end{equation}
We now set $z_{k}=u_{n_{k}}-q_{k}-w_{k}$. Note that in particular $|z_{k}| \leq |u_{n_{k}}|.$ We have 
$$\int |z_{k}|^{2} \leq \int_{\frac{t_{k}}{2}\leq |x-y(u_{n_{k}},\frac{t_{k}}{2})|\leq t_{k}} |u_{n_{k}}|^{2}=\int_{|x-y(u_{n_{k}},\frac{t_{k}}{2})|
\leq t_{k}} |u_{n_{k}}|^{2}-\int_{|x-y(u_{n_{k}},\frac{t_{k}}{2})|\leq \frac{t_{k}}{2}} |u_{n_{k}}|^{2}$$
$$\leq \int_{|x-y(u_{n_{k}},t_{k})| \leq t_{k}} |u_{n_{k}}|^{2}-\int_{|x-y(u_{n_{k}},\frac{t_{k}}{2})|\leq \frac{t_{k}}{2}} |u_{n_{k}}|^{2}=
\rho(u_{n_{k}},t_{k})-\rho(u_{n_{k}},\frac{t_{k}}{2})$$
and again by Lemma 6 we have
\begin{equation}
\lim_{k \rightarrow \infty}\|z_{k}\|_{L^{2}}^{2}=0.
\end{equation}
By Cauchy-Schwartz inequality and the above we have that
$$\lim_{k \rightarrow \infty}\int u_{n_{k}}\bar{z}_{k}=0$$
But now by $(24)$, $(30)$, $(31)$ and some trivial algebra we get after integration that
$$\lim_{k \rightarrow \infty}\|w_{k}\|_{L^{2}}^{2}=b-\mu$$
and (27) follows. Also note that $z_{k}$ is bounded in $H^{1}$ and converges to $0$ in $L^{2}$, 
and by Gagliardo-Nirenberg inequality, in $L^{p}$ for any $2 \leq p<\infty$. Moreover one can easily verifies that
$$\Bigr||u_{n_{k}}|^{p}-|q_{k}|^{p}-|w_{k}|^{p}\Bigr| \leq C |u_{n_{k}}|^{p-1}|z_{k}|$$
and by Cauchy-Schwartz since $z_{k}$ tends to $0$ in $L^{p}$, $(29)$ follows. Finally $(28)$ follows easily from the initial assumptions, 
Cauchy-Schwartz inequality and the easy calculation
$$|\nabla u_{n_{k}}|^{2}-|\nabla q_{k}|^{2}-|\nabla w_{k}|^{2}=|\nabla u_{n_{k}}|^{2}(1-\theta_{k}^{2}-\phi_{k}^{2})-|u_{n_{k}}|^{2}(|\nabla \theta_{k}|^{2}
+|\nabla \phi_{k}|^{2})$$
$$-Re(\bar{u}_{n_{k}}\nabla u_{n_{k}})\cdot \nabla (\theta_{k}^{2}+\phi_{k}^{2}) \geq -\frac{C}{t_{k}^{2}}|u_{n_{k}}|^{2}-\frac{C}{t_{k}}|u_{n_{k}}|\ |\nabla u_{n_{k}}|.$$
\end{proof}
{\bf Proof of Theorem 3.}
\begin{proof}
Define the blowup parameters:
\\
$$\lambda(t)=\|u(t)\|_{H^{s}},\ \ \Lambda(t)=\sup_{0\leq \tau \leq t}\lambda(\tau)$$
\\
$$\sigma(t)=\|I_{N}\langle \nabla\rangle u(t)\|_{L^{2}},\ \ \Sigma(t)=\sup_{0\leq \tau \leq t}\sigma(\tau)$$
\\
Let $\{t_{n}\}_{n=1}^{\infty}$ be a sequence such that $t_{n}\uparrow T^{\star}$ and for each $t_{n}$ we have
\\
$$\|u(t_{n})\|_{H^{s}}=\Lambda(t_{n})$$
and with $u(t_{n})=u_{n}$ we define 
$$I_{N}u_{n}=I_{N(t_{n})}u(t_{n}).$$
We rescale these as follows
$$v_{n}(x)=\frac{1}{\sqrt{\sigma_{n}}}I_{N}u_{n}(\frac{x}{\sigma_{n}})$$
where
$$\sigma_{n}=\|I_{N}\langle \nabla \rangle u_{n}\|_{L^{2}}=\sigma(t_{n})$$
\\
Note that for these sequences(let's call them maximizing) we have that
$$\Lambda(t_{n})\leq \sigma_{n}$$
where $\sigma_{n}\rightarrow \infty$ as $n\rightarrow \infty$. It is important to note that we are in the
 blow-up regime and thus 
$$\|u_{0}\|_{L^{2}} \geq \|Q\|_{L^{2}}.$$ 
Moreover the $L^{2}$ part of $v_{n}$ is bounded uniformly in $n$. This is because
$$\|v_{n}\|_{L^{2}}=\|I_{N}u_{n}\|_{L^{2}} \leq \|u_{n}\|_{L^{2}}=\|u(t_{n})\|_{L^{2}}=\|u_{0}\|_{L^{2}}.$$
Thus in the limit as $n \rightarrow \infty$ we have
$$\lim_{n \rightarrow \infty}\|\nabla v_{n}\|_{L^{2}}=1$$
Also since $N(t_{n})$ goes to infinity as $n \rightarrow \infty$
$$ \lim_{n \rightarrow \infty}\|v_{n}\|_{L^{2}}=\lim_{n \rightarrow \infty}\|I_{N}u_{n}\|_{L^{2}}=\lim_{n \rightarrow \infty}\|u(t_{n})\|_{L^{2}}=
\|u_{0}\|_{L^{2}} \geq \|Q\|_{L^{2}}.$$
In addition by Proposition 3 we have that
$$|E(v_{n})|=\frac{1}{\sigma_{n}^{2}}|E(I_{N}u_{n})|\leq C\sigma_{n}^{-2}\Lambda^{p(s)}(t_{n})\leq C\Lambda^{p(s)-2}(t_{n})$$
and thus 
$$\lim_{n\rightarrow \infty}E(v_{n})=0.$$
since $p(s)<2$. This allows another way to prove that 
$$\lim_{n \rightarrow \infty}\|v_{n}\|_{L^{2}}\geq \|Q\|_{L^{2}}$$ since by the optimality of Gagliardo-Nirenberg inequality we have
$$E(v_{n}) \geq \frac{1}{2}(1-\frac{\|v_{n}\|_{L^{2}}^{4}}{\|Q\|_{L^{2}}^{4}})\|\nabla v_{n}\|_{L^{2}}^{2}$$ and in the limit as $n \rightarrow \infty$
 we get 
$$\lim_{n \rightarrow \infty}\|v_{n}\|_{L^{2}}\geq \|Q\|_{L^{2}}.$$
We collect the three important relations that we have
\begin{equation}
\lim_{n \rightarrow \infty}\|v_{n}\|_{L^{2}}=\|u_{0}\|_{L^{2}}\geq \|Q\|_{L^{2}},
\end{equation}
\begin{equation}
\lim_{n \rightarrow \infty}\|\nabla v_{n}\|_{L^{2}}=1,
\end{equation}
\begin{equation}
\lim_{n \rightarrow \infty}E(v_{n})=0.
\end{equation}
With the help of (32), (33), and (34) we will conclude that
\\
\\
{\bf Claim:}
$$\mu(\{v_{n}\}_{n \geq 0}) \geq \|Q\|_{L^{2}}^{2}.$$ 
First assuming the claim and revisiting the statement of the theorem it is enough to prove that
 for any $\epsilon>0$ we have
\\
$$\lim_{n\rightarrow \infty} \|u(t_{n})\|_{L_{\{|x-z_{n}|<(T^{\star}-t_{n})^{\frac{s}{2}}\gamma(T^{\star}-t_{n})\}}^{2}} \geq \|Q\|_{L^{2}}-\epsilon$$
Note that since $N(t_{n})$ goes to $\infty$ we have
\\ 
$$\lim_{n\rightarrow \infty}\|u(t_{n})\|_{L^{2}}=
\lim_{n\rightarrow \infty}\|I_{N(t_{n})}u(t_{n})\|_{L^{2}}$$
Now given $\epsilon>0$, the relation
$$\mu(\{v_{n}\}_{n\geq 0}) \geq \|Q\|_{L^{2}}^{2}$$
\\
by Lemma 6 implies that there exist a $T \in {\bf R}$ such that
\\
$$\rho(v_{n},T)\geq \|Q\|_{L^{2}}^{2}-\epsilon$$
\\
for large $n$. Note that $\rho$ is a nondecreasing function of $t$ and it is crucial to find a fixed $T$ 
such that $\rho(v_{n},T)\geq \|Q\|_{L^{2}}^{2}-\epsilon$ holds for any large $n$, for the given $\epsilon$. That this $T$ is independent of all the $n's$ after some $n$ large,
up to a subsequence is guaranteed by the two parts of Lemma 6. Thus the same $T$ works for all large $n$.
\\
\\
Now setting $y_{n}=y(v_{n},T)$ defined by Lemma 4 we have that
\\
$$\liminf_{n\rightarrow \infty}\|v(t_{n})\|_{L_{\{|x-y_{n}|<T\}}^{2}}\geq \|Q\|_{L^{2}}^{2}-\epsilon$$
\\
and up to a subsequence
\\
\begin{equation}
\lim_{n\rightarrow \infty}\|v(t_{n})\|_{L_{\{|x-y_{n}|<T\}}^{2}}\geq \|Q\|_{L^{2}}^{2}-\epsilon
\end{equation}
But 
$$\lim_{n\rightarrow \infty}\|v(t_{n})\|_{L_{\{|x-y_{n}|<T\}}^{2}}=
\lim_{n\rightarrow \infty}\|\frac{1}{\sqrt{\sigma_{n}}}I_{N}u_{n}(\frac{x}{\sigma_{n}})\|_{L_{\{|x-y_{n}|<T\}}^{2}}=$$
\\
$$\lim_{n\rightarrow \infty}\|I_{n}u_{n}(x)\|_{L_{\{|x-\frac{y_{n}}{\sigma_{n}}|<\frac{T}{\sigma_{n}}\}}^{2}}=
\lim_{n\rightarrow \infty}\|u_{n}(x)\|_{L_{\{|x-\frac{y_{n}}{\sigma_{n}}|<\frac{T}{\sigma_{n}}\}}^{2}}=
\lim_{n\rightarrow \infty}\|u_{n}(x)\|_{L_{\{|x-z_{n}|<\frac{T}{\sigma_{n}}\}}^{2}}$$
\\
where $z_{n}=\frac{y_{n}}{\sigma_{n}}$. Moreover note that $\frac{T}{\sigma_{n}}\rightarrow 0$
and that $\sigma_{n}$ goes to infinity at least as fast as $(T^{\star}-t)^{-\frac{s}{2}}$. Thus there exists a funtion $\gamma(x) \uparrow \infty$ as 
$x\downarrow 0$ such that
\\
$$\lim_{n\rightarrow \infty} \|u(t_{n})\|_{L_{\{|x-z_{n}|<(T^{\star}-t_{n})^{\frac{s}{2}}\gamma(T^{\star}-t_{n})\}}^{2}}\geq
\lim_{n\rightarrow \infty}\|u_{n}(x)\|_{L_{\{|x-z_{n}|<\frac{T}{\sigma_{n}}\}}^{2}}$$
\\
$$=\lim_{n\rightarrow \infty}\|v(t_{n})\|_{L_{\{|x-y_{n}|<T\}}^{2}} \geq \|Q\|_{L^{2}}^{2}-\epsilon$$ 
\\
by equation $(35)$ and the proof is complete.
\\
\\
{\bf Proof of Claim}
\\
\\
We prove the claim by contradiction. We claim that there exists $\delta>0$ with the following property. If $(v_{n})_{n\geq 0} \in H^{1}$ is such that
\begin{equation}
\lim_{n \rightarrow \infty}\|v_{n}\|_{L^{2}}^{2}=\|u_{0}\|_{L^{2}}^{2}
\end{equation}
\begin{equation}
0<\liminf_{n \rightarrow \infty}\|\nabla v_{n}\|_{L^{2}}\leq \limsup_{n \rightarrow \infty}\|\nabla v_{n}\|_{L^{2}}<\infty
\end{equation}
\begin{equation}
\limsup_{n \rightarrow \infty}E(v_{n})\leq 0
\end{equation}
and
\begin{equation}
\mu(\{v_{n}\}_{n \geq 0})<\|Q\|_{L^{2}}^{2}
\end{equation}
then there exists a sequence $(\tilde{v}_{n})_{n\geq 0} \in H^{1}$ satisfying (37),(38),(39) and such that 
$$\lim_{n \rightarrow \infty}\|\tilde{v}_{n}\|_{L^{2}}^{2}=\|u_{0}\|_{L^{2}}^{2}-\beta$$ 
for some $\beta > \delta$. Clearly the sequence $v_{n}$ of Theorem 3 satisfies (36), (37), and (38). But then by Galiardo-Nirenberg and
(37),(38),(39) we have that
$$\mu(\{v_{n}\}_{n \geq 0})<\|u_{0}\|_{L^{2}}^{2}-\delta $$
If we apply the above proceedure $k-$times we get $\mu(\{v_{n}\}_{n \geq 0})<\|u_{0}\|_{L^{2}}^{2}-k \delta $ which for large $k$ is absurd. 
Thus it suffices to prove the claim. We apply Lemmas 6 and 7 to the sequence $(v_{n})_{n\geq 0}$ and we consider the corresponding sequences
$(q_{n})_{n\geq 0}$ and $(w_{n})_{n\geq 0}$.We set 
$$\delta=(\frac{3}{K})^{\frac{1}{2}}>0$$ 
where $K$ is given in Lemma 5. We first show that
$$\mu(\{v_{n}\}_{n \geq 0})\geq \delta.$$
By Lemma 5 and the definition of $\delta$ we have that
$$E(v_{n_{k}})\geq \frac{1}{2}\left( 1-(\frac{\rho(v_{n_{k}},t_{k})}{\delta})^{2}\right)\int |\nabla v_{n_{k}}|^{2}-\frac{K}{6t_{k}^{2}}\rho(v_{n_{k}},t_{k})^{2}$$
Now if we assume by contradiction that $\mu<\delta$, then we obtain by letting $k \rightarrow \infty$, applying the second part of Lemma 6, and (37) that
up to a subsequence
$$\limsup_{n \rightarrow \infty}E(v_{n})\geq \frac{1}{2}\left( 1-(\frac{\mu}{\delta})^{2}\right)\liminf_{n \rightarrow \infty} \int |\nabla v_{n}|^{2}>0$$
which is absurd. Now since by $(25)$ we know that $|w_{k}| \leq |v_{n_{k}}|$ we have by the second part of Lemma 6 and (39)
$$\mu((w_{k})_{k\geq 0}) \leq \mu <\|Q\|_{L^{2}}^{2}.$$
This proves that  $(w_{k})_{k\geq 0}$ satisfies (39). 
Also by $(27)$, $(39)$ and Gagliardo-Nirenberg inequality we know that there exists a $\sigma>0$ such that for $k$ large
\begin{equation}
E(q_{k})\geq \sigma \|\nabla q_{k}\|_{L^{2}}^{2}
\end{equation}
On the other hand by (28) and (29) we have that
\begin{equation}
\liminf_{k \rightarrow \infty}\{E(v_{n_{k}})-E(q_{k})-E(w_{k}))\}\geq 0
\end{equation}
and thus 
$$\limsup_{k \rightarrow \infty}E(w_{k}) \leq 0.$$
This proves that  $(w_{k})_{k\geq 0}$ satisfies (38). By (26) and (37) we easily get that
$$\|\nabla w_{k}\|_{L^{2}} \leq C\|v_{n_{k}}\|_{H^{1}}< \infty$$
Finally we show the last property (37), namely that
$$\liminf_{k \rightarrow \infty}\|\nabla w_{k}\|_{L^{2}}>0$$
We argue again by contradiction and assume that there exists a sequence which we still denote by $(w_{k})_{k \geq 0}$ such that 
$\lim_{k \rightarrow \infty}\|\nabla w_{k}\|_{L^{2}}=0$. But then trivially follows that $E(w_{k}) \rightarrow 0$ as $k \rightarrow \infty$ 
and thus by (38), (40), and (41) we get that
$$\lim_{k \rightarrow \infty}\|\nabla q_{k}\|_{L^{2}}=0$$
But then by (29) and the fact that $\|w_{k}\|_{L^{6}},\|q_{k}\|_{L^{6}} \rightarrow 0$ we deduce that $\lim_{k \rightarrow \infty}\|v_{n_{k}}\|_{L^{6}}=0$
and thus
$$\limsup_{k \rightarrow \infty}E(v_{n_{k}})>0$$
which contradicts (38). Now setting
$$\tilde{v}_{k}=\frac{\sqrt{\|u_{0}\|_{L^{2}}^{2}-\mu}}{\|w_{k}\|_{L^{2}}}w_{k}$$ 
we see that with the help of (27) the sequence $(\tilde{v}_{n})_{n \geq 0}$ satisfies (37), (38), (39) and that
$$\lim_{k \rightarrow \infty}\|\tilde{v}_{k}\|_{L^{2}}^{2}=\|u_{0}\|_{L^{2}}^{2}-\mu \leq \|u_{0}\|_{L^{2}}^{2}-\delta$$
and we are done.
\end{proof}
{\bf Remark.} In dimensions $n\geq 2$ with the additional assumption of radial symmetry on the initial data, the solution of the equivalent $L^{2}$ critical
 Schr\"odinger equation, satisfies the conclusion of Theorem 3 with $z(t)\equiv 0$.
\\
\\
We define the ``modified energy'' for the system $(16)$ as
$$E(Iu)(t)=\frac{1}{2}\int{|Iu_{x}(t)|^{2}dx}-\frac{1}{6} \int{|Iu(t)|^{6}dx}$$
This ``energy'' functional is not conserved but we can show that its time derivative decays with respect to a large parameter $N$. The next proposition quantifies
 the increament of this functional on $[0,\delta]$.
\begin{proposition}
Let $u$ be an $H^{1}$ solution of $(2).$ Then
\\
$$E(Iu)(\delta)-E(Iu)(0)=$$
$$Im\left(\int_{0}^{\delta} \int I\bar{u}_{xx}\left(I(|u|^{4}u)-Iu|Iu|^{4}\right)dxdt\right)+
Im\left(\int_{0}^{\delta}\int I(|u|^{4}u)\left(I(|u|^{4}u)-Iu|Iu|^{4}\right)dxdt\right).$$
\end{proposition}
\begin{proof}
The derivative of the ``modified energy'' is:
\\
$$\frac{dE}{dt}(Iu)=Im\left(\int I\bar{u}_{xx}\left(I(|u|^{4}u)-Iu|Iu|^{4}\right)dx\right)+Im\left(\int I(|u|^{4}u)\left(I(|u|^{4}u)-Iu|Iu|^{4}\right)dx \right).$$
\\
But then Proposition 4 follows immediately by applying the fundamental theorem of calculus.
\end{proof}
By the previous formal identity we can deduce the desired decay of the ``modified energy''.
\begin{proposition}
For any Schwartz function $u$ we have that
$$E(Iu)(\delta)-E(Iu)(0) \lesssim 
\delta^{\frac{1}{4}-}N^{-\frac{3}{2}+}\|Iu\|_{X_{\delta}^{1,1/2+}}^{6}+\delta^{\frac{1}{2}-}N^{-2+}\|Iu\|_{X_{\delta}^{1,1/2+}}^{10}$$
\end{proposition}
\begin{proof}
First we establish
$$|\int_{0}^{\delta} \int I\bar{u}_{xx}\left(I(|u|^{4}u)-Iu|Iu|^{4}\right)dxdt| \lesssim \delta^{\frac{1}{4}-}N^{-\frac{3}{2}+}\|Iu\|_{X_{\delta}^{1,1/2+}}^{6}$$
\\
or by Plancherel's theorem that
\\
$$|\int_{0}^{\delta}\int_{\Gamma_{6}}\frac{\xi_{1}^{2}}{<\xi_{1}>}
(\frac{m(\xi_{2}+...+\xi_{6})-m(\xi_{2})...m(\xi_{6})}{m(\xi_{2})...m(\xi_{6})})
\hat{u}(\xi_{1},t)...\hat{\bar{u}}(\xi_{6},t)d\xi dt|$$
\begin{equation}
\lesssim \delta^{1/4-}N^{-\frac{3}{2}+}\|u\|_{X^{1,1/2+}}^{5}\|u_{1}\|_{X^{0,1/2+}}
\end{equation}
where $\Gamma_{6}$ denotes the hyperplane $\xi_{1}+\xi_{2}+...+\xi_{6}=0$, and $u_{1}$ is the function that corresponds 
on the Fourier side to the frequency $\xi_{1}$.
\\
\\
$\bf{Remarks.}$
\\
\\
$\bf{1.}$
Let us denote $N_{i}\sim |\xi_{i}|$ and $N_{max}\sim |\xi|_{max}$,
$N_{med}\sim |\xi|_{med}$ where  $|\xi|_{max}$, $|\xi|_{med}$ is the
largest and second largest of the $|\xi_{i}|$. 
If all $|\xi_{i}| \ll N$ then the parenthesis above is zero and there is nothing to prove.
Thus since the $\xi_{i}$ are related by
 $\xi_{1}+\xi_{2}+...+\xi_{6}=0$ we have that $|\xi|_{max} \sim |\xi|_{med} \gtrsim N.$ We also write $m_{i}$ for $m(\xi_{i})$ and $m_{ij}$ for $m(\xi_{i}+\xi_{j})$.
\\
\\
$\bf{2.}$ Our strategy from now on is to break all the functions into a sum of dyadic constituents $\psi_{j}$, each
with frequency support $<\xi> \sim 2^{j}, j=0,...$
Then we pull the absolute value of the symbols out of the integral, estimating it pointwise. After
bounding the multiplier, 
the remaining integrals involving the pieces $\psi_{j}$ are estimated by reversing the Plancherel formula and
using duality, H\"older's inequality and Strichartz's estimates. We can sum over all the frequency pieces
$\psi_{j}$ as long as we keep always a factor $N_{max}^{-\epsilon}$ inside the summation.
\\
\\
$\bf{3.}$ Since in all of the estimates that we establish from now on,
the right hand side is in terms of the $X^{s,b}$ norms and the
$X^{s,b}$ spaces depend only on the absolute value of the Fourier
transform, we can assume without loss of generality that the Fourier
transform of all the functions in the estimates are real and positive.
\\
\\
$\bf{4. }$ Note also that $\frac{N_{1}^{2}}{\langle N_{1} \rangle} \leq N_{1}$.
\\
\\
Since our analysis as we mentioned before do not rely upon the complex conjugate structure of the left hand side, there is a symmetry under
 the interchange of the indices and thus we can assume that
$$N_{2} \geq N_{3}\geq ... \geq N_{6}$$
Case 1: Let $N \gg N_{2}$. Then 
\\
$$\frac{m(\xi_{2}+...+\xi_{6})-m(\xi_{2})...m(\xi_{6})}{m(\xi_{2})...m(\xi_{6})}=0$$
and there is nothing to prove.
\\
\\
Case 2: $N_{2} \gtrsim N \gg N_{3} \geq ...\geq N_{6}$. This forces $N_{1}\sim N_{2}$ on $\Gamma_{6}$. 
\\
\\
But then by the mean value theorem we have
\\
$$\Bigr| \frac{m(\xi_{2}+...+\xi_{6})-m(\xi_{2})...m(\xi_{6})}{m(\xi_{2})...m(\xi_{6})}\Bigr|=
\Bigr|\frac{m(\xi_{2})-m(\xi_{1})}{m(\xi_{2})}\Bigr| \lesssim$$
\\
$$\Bigr|\frac{\nabla m(\xi_{2})\cdot (\xi_{3}+...+\xi_{6})}{m(\xi_{2})}\Bigr| \lesssim \frac{N_{3}}{N_{2}}$$
\\
Now by undoing Plancherel's theorem, using Cauchy-Schwartz inequality, apply the Strichartz estimates and using Lemmas 2 and 3 we have that the left hand side of 
$(42)$ is
$$\lesssim \frac{N_{1}N_{3}}{N_{2}N_{1}^{1/2}}\|(D^{1/2}u_{1})u_{3}\|_{L_{t}^{2}L_{x}^{2}}\|u_{2}u_{4}u_{5}u_{6}\|_{L_{t}^{2}L_{x}^{2}}$$
$$\lesssim \frac{N_{3}}{N_{1}^{1/2}}\|u_{1}\|_{X_{\delta}^{0,1/2+}}\|u_{3}\|_{X_{\delta}^{0,1/2+}}
\prod_{j=4}^{6}\|u_{j}\|_{L_{t}^{\infty}L_{x}^{\infty}}\|u_{2}\|_{L_{t}^{2}L_{x}^{2}}$$
$$\lesssim \delta^{1/2-}\frac{N_{3}}{N_{1}^{1/2}}\prod_{j=1}^{3}\|u_{j}\|_{X_{\delta}^{0,1/2+}}\prod_{j=4}^{6}\|u_{j}\|_{X_{\delta}^{1/2,1/2+}}$$
where in the last inequality we also used equation $(9)$ in its dyadic form. Comparing with $(42)$ we see that it is enough to have
\\
$$\delta^{1/2-}\frac{N_{3}(N_{4}N_{5}N_{6})^{1/2}}{N_{1}^{1/2}}\lesssim \delta^{1/4-}N^{-\frac{3}{2}+}N_{max}^{-\epsilon}\langle N_{2} \rangle...\langle N_{6} \rangle$$
which is true. Note that in the process we summed the Littlewood-Paley pieces, using the factor $N_{max}^{-\epsilon}$.
\\
\\
Case 3: $N_{2}\geq N_{3}\gtrsim N$. 
\\
\\
In this case we use the crude estimate
\\
$$|1-\frac{m_{1}}{m_{2}...m_{6}}| \lesssim \frac{m_{1}}{m_{2}...m_{6}}$$ 
\\
Since it is impossible to have 
$N_{1} \gg N_{med}=N_{2}$, we can divide this case into two subcases.
\\
\\
a) $N_{1}\sim N_{2}\geq N_{3}\gtrsim N$. Now we are starting comparing the different frequencies in order to be able to apply Lemma 3.
 Note that $m_{1}\sim m_{2}$.
\\
\\
i) Suppose first that $N_{2}\gg N_{3}$. Without loss of generality we can assume that $N_{4} \leq N$. This is because in case that one
 of the $N_{4}, N_{5},N_{5}$ are $\gtrsim N$ the estimate is even easier and the decay is greater. For the suspicious reader that might object the
previous argument because
 of the presence of $m_{4}m_{5}m_{6}$ in the denominator we comment that for $N_{j} \gtrsim N$ we have that 
$$\frac{1}{m_{j}N_{j}^{1/2}} \lesssim \frac{1}{N^{1/2}}$$ and indeed we can get a better decay. From now on
 we will use this heuristic without any comment. Thus we can apply Cauchy-Schwartz  and Lemma 3
 and the left hand side of $(42)$ is
\\
$$\lesssim \frac {N_{1}}{m_{3}N_{1}^{1/2}}\|(D^{1/2}u_{1})u_{3}\|_{L_{t}^{2}L_{x}^{2}}\|u_{2}\|_{L_{t}^{2}L_{x}^{2}}
\prod_{j=4}^{6}\|u_{j}\|_{L_{t}^{\infty}L_{x}^{\infty}}$$ 
$$\lesssim \delta^{1/2-}\frac{N_{1}N_{3}^{1/2}}{m_{3}N_{3}^{1/2}N_{1}^{1/2}}\prod_{j=1}^{3}\|u_{j}\|_{X_{\delta}^{0,1/2+}} \cdot 
\prod_{j=4}^{6}\|u_{j}\|_{X_{\delta}^{1/2,1/2+}}$$
$$\lesssim \delta^{1/2-}\frac{N_{1}N_{3}^{1/2}}{N^{1/2}N_{1}^{1/2}}\prod_{j=1}^{3}\|u_{j}\|_{X_{\delta}^{0,1/2+}} \cdot \prod_{j=4}^{6}\|u_{j}\|_{X_{\delta}^{1/2,1/2+}}.$$
\\
Comparing with $(42)$ we see that it is enough to have
$$\delta^{1/2-}\frac{N_{1}N_{3}^{1/2}}{N^{1/2}N_{1}^{1/2}} \lesssim \delta^{1/4-}N^{-\frac{3}{2}+}
N_{max}^{-\epsilon}(\langle N_{4}\rangle...\langle N_{6} \rangle)^{1/2}\langle N_{2} \rangle \langle N_{3} \rangle$$
\\
which is true.
\\
\\
ii) Now assume that $N_{2} \sim N_{3}$ and by the comment in case i) the worst case is when $N_{5},N_{6} \leq N$ which we assume without loss of generality.
. In this case we compare $N_{3}$ with $N_{4}$. In case that $N_{3} \sim N_{4}$ the estimate
is easy since
$$\frac{N_{1}m_{1}}{m_{2}...m_{6}} \lesssim \frac{N_{1}N_{3}^{1/2}}{m_{3}m_{4}N_{3}^{1/2}} \lesssim \frac{N_{1}N_{3}^{1/2}}{N^{1/2}}$$
\\
and thus the left hand side of $(42)$ is
$$\lesssim \frac{N_{1}N_{3}^{1/2}}{N^{1/2}} \prod_{j=1}^{4}\|u_{j}\|_{L_{t}^{4}L_{x}^{4}}\cdot 
\|u_{5}\|_{L_{t}^{\infty}L_{x}^{\infty}}\|u_{6}\|_{L_{t}^{\infty}L_{x}^{\infty}}$$
$$\lesssim \delta^{1/2-}\frac{N_{1}N_{3}^{1/2}}{N^{1/2}} 
\prod_{j=1}^{4}\|u_{j}\|_{X_{\delta}^{0,1/2+}} \cdot \|u_{5}\|_{X_{\delta}^{1/2,1/2+}}\|u_{6}\|_{X_{\delta}^{1/2,1/2+}}.$$
Comparing with $(42)$ it is enough to have
\\
$$\delta^{1/2-}\frac{N_{1}N_{3}^{1/2}(N_{5}N_{6})^{1/2}}{N^{1/2}} \lesssim \delta^{1/4-}N^{-\frac{3}{2}+}N_{max}^{-\epsilon}\langle N_{2}\rangle...\langle N_{6} \rangle$$
\\
which is true. If $N_{3} \gg N_{4}$ again without loss of generality we assume that $N_{4} \leq N$ and we apply Lemma 3. Moreover
\\
$$\Bigr|\frac{N_{1}m_{1}}{m_{2}...m_{6}}\Bigr| \lesssim \frac{N_{1}N_{3}^{1/2}}{m_{3}N_{3}^{1/2}} \lesssim \frac{N_{1}N_{3}^{1/2}}{N^{1/2}}$$
and thus the left hand side of $(42)$ is
$$\lesssim \frac{N_{1}N_{3}^{1/2}}{N^{1/2}N_{1}^{1/2}}
\|(D^{1/2}u_{1})u_{4}\|_{L_{t}^{2}L_{x}^{2}}\|u_{2}u_{3}\|_{L_{t}^{2}L_{x}^{2}}\|u_{5}\|_{L_{t}^{\infty}L_{x}^{\infty}}\|u_{6}\|_{L_{t}^{\infty}L_{x}^{\infty}}$$
$$\lesssim \frac{N_{1}N_{3}^{1/2}}{N^{1/2}N_{1}^{1/2}}
\|u_{1}\|_{X_{\delta}^{0,1/2+}}\|u_{4}\|_{X_{\delta}^{0,1/2+}}\|u_{2}\|_{L_{t}^{6}L_{x}^{6}}\|u_{3}\|_{L_{t}^{3}L_{x}^{3}}
\|u_{5}\|_{X_{\delta}^{1/2,1/2+}}\|u_{6}\|_{X_{\delta}^{1/2,1/2+}}$$
$$\lesssim \delta^{1/4-}\frac{N_{1}N_{3}^{1/2}}{N^{1/2}N_{1}^{1/2}} 
\prod_{j=1}^{4}\|u_{j}\|_{X_{\delta}^{0,1/2+}}\|u_{5}\|_{X_{\delta}^{1/2,1/2+}}\|u_{6}\|_{X_{\delta}^{1/2,1/2+}}.$$
\\
where we used Lemmas 2 and 3 and equations $(7)$, $(9)$ and $(11)$. Comparing with $(42)$ it is enough to have
\\
$$\delta^{1/4-}\frac{N_{1}N_{3}^{1/2}(N_{5}N_{6})^{1/2}}{N^{1/2}N_{1}^{1/2}} \lesssim \delta^{1/4-}N^{-\frac{3}{2}+}N_{max}^{-\epsilon}\langle N_{2} \rangle...
\langle N_{6} \rangle$$
which is true.
\\
\\
b)\ $N_{2} \sim N_{3} \gtrsim N$ and $N_{2} \gg N_{1}$. Since $N_{1}$ is in the numerator on the left hand side of equation $(20)$ this case
 is easier than the previous and similar analysis gives the same (or ever better) bounds as in a). The details are omitted. 
\\
\\
To conclude the proof of Proposition 5 it remains to show that
\\
$$\Bigr|\int_{0}^{\delta}\int I(|u|^{4}u)\left(I(|u|^{4}u)-Iu|Iu|^{4}\right)dxdt\Bigr| \lesssim \delta^{\frac{1}{2}-}N^{-2+}\|Iu\|_{X_{\delta}^{1,1/2+}}^{10}$$
\\
By Plancerel's theorem
\begin{equation}
|\int_{0}^{\delta}\int_{\Gamma_{10}}m_{12345}\{m_{678910}-m_{6}...m_{10}\}\hat{u}(\xi_{1},t)...\hat{\bar{u}}(\xi_{10},t)d\xi dt| 
\lesssim \delta^{\frac{1}{2}}N^{-2+}\|Iu\|_{X_{\delta}^{1,1/2+}}^{10}
\end{equation}
As we noted before if $N_{max} \ll N$ the multiplier is zero so we assume that 
$$N_{max} \sim N_{med} \gtrsim N.$$
In addition since $m(\xi) \leq 1$ we have that
\\
$$|m_{12345}\{m_{678910}-m_{6}...m_{10}\}| \lesssim C.$$
\\
Finally the last pointwise estimate that we use is the following
\\
$$\frac{1}{m_{max}N_{max}} \lesssim N^{-1}$$
\\
which follows easily since $N_{max} \gtrsim N$. The left hand side of $(21)$ is
\\ 
$$\lesssim \int_{0}^{\delta}\int_{\Gamma_{10}}\prod_{j=1}^{10}\hat{u}(\xi_{j},t)d\xi dt 
\lesssim \int_{0}^{\delta}\int_{\Gamma_{10}}\frac{m_{max}N_{max}\hat{u}_{max}\cdot 
m_{med}N_{med}\hat{u}_{med}}{m_{max}N_{max}m_{med}N_{med}}\prod_{j\neq j_{max}, j_{med}}\hat{u}(\xi_{j},t)d\xi dt$$
$$\lesssim N^{-2+}N_{max}^{-\epsilon}\int_{0}^{\delta}\int_{\Gamma_{10}}\widehat{DIu_{max}}\cdot \widehat{DIu_{med}}\prod_{j\neq j_{max}, j_{med}}\hat{u}(\xi_{j},t)d\xi dt.$$
Now reversing Prancherel's Theorem and use the the following estimates
\\
$$\|u\|_{L_{t}^{6}L_{x}^{6}} \lesssim \|u\|_{X_{\delta}^{0,1/2+}}$$
\\
$$\|u\|_{L_{t}^{3}L_{x}^{3}} \lesssim \delta^{1/4-}\|u\|_{X_{\delta}^{0,1/2+}}$$
\\
$$\|u\|_{X_{\delta}^{0,1/2+}}\lesssim \|Iu\|_{X_{\delta}^{1,1/2+}}$$
\\
$$\|u\|_{X_{\delta}^{1/2,1/2+}}\lesssim \|Iu\|_{X_{\delta}^{1,1/2+}}$$
\\
we get that the left hand side of $(43)$ is
\\
$$\lesssim  N^{-2+}N_{max}^{-\epsilon}\|JIu_{max}\|_{L_{t}^{6}L_{x}^{6}} \cdot \|JIu_{med}\|_{L_{t}^{6}L_{x}^{6}} \|u\|_{L_{t}^{3}L_{x}^{3}}^{2}
\|u\|_{L_{t}^{\infty}L_{x}^{\infty}}^{6}\lesssim$$
\\ 
$$\delta^{1/2-}N^{-2+}N_{max}^{-\epsilon}\|JIu\|_{X_{\delta}^{0,1/2}}^{2}\|u\|_{X_{\delta}^{0,1/2}}^{2}\|u\|_{X_{\delta}^{1/2,1/2}}^{6} 
\lesssim \delta^{\frac{1}{2}-}N^{-2+}\|Iu\|_{X_{\delta}^{1,1/2+}}^{10}$$
where in the process we sum the different Littlewood-Paley pieces, taking advantage of the factor $N_{max}^{-\epsilon}$.
\end{proof}
Now we are finally ready to prove Proposition 3.
\begin{proof}
When $s=1$ we can choose $N(T)=+\infty$ and thus $I_{N(T)}=1$ and the proposition is true with $p(s)=0$ since the energy 
is conserved and the kinetic energy blows up as time approaches $T^{\star}$. Therefore we can fix $10/11<s<1$ and take $T$ near $T^{\star}$. Now let
 $N=N(T)$ to be chosen later in the argument. Recall that $\delta \sim (\Sigma(T))^{-4-\epsilon}$ gives the time of the local well-posedness. Thus
 if we divide the interval $[0,T]$ into $\frac{T}{\delta}$-subintervals of size $\sim \delta$, the local well-posedness result uniformly
 applies. Moreover for any $t$ in this subinterval we have that
\\
$$\|I\langle \nabla\rangle u(t)\|_{L^{2}}=\sigma(t) \leq \Sigma(T)$$
\\
The next step is to iterate the almost conservation of the energy. It is apparent that after $\frac{T}{\delta}$-steps the growth of the
 modified energy is
$$E(Iu(T)) \lesssim E(Iu(0))+\frac{T^{\star}}{\delta}\{\delta^{\frac{1}{4}-}N^{-\frac{3}{2}+}\Sigma(T)^{6}+\delta^{\frac{1}{2}-}N^{-2+}\Sigma(T)^{10}\}$$
\\
$$\leq N^{2(1-s)}\lambda(0)+\frac{T^{\star}}{\delta}\{\delta^{\frac{1}{4}-}N^{-\frac{3}{2}+}\Sigma(T)^{6}+\delta^{\frac{1}{2}-}N^{-2+}\Sigma(T)^{10}\}$$
\\
$$\lesssim N^{2(1-s)}+\{\delta^{-\frac{3}{4}-}N^{-\frac{3}{2}+}\Sigma(T)^{6}+\delta^{-\frac{1}{2}-}N^{-2+}\Sigma(T)^{10}\}$$
\\ 
$$\lesssim N^{2(1-s)}+N^{-\frac{3}{2}+}\Sigma(T)^{9+}+N^{-2+}\Sigma(T)^{12+}$$
\\
where in the third inequality we dismiss the irrelevant constants. Now if we switch from $\Sigma(T)$ to $\Lambda(T)$ we have
\\ 
$$E(Iu(T)) \lesssim N^{2(1-s)}+N^{-\frac{3}{2}+}N^{9(1-s)+}\Lambda(T)^{9+}+N^{-2+}N^{12(1-s)+}\Lambda(T)^{12+}.$$
\\
We know choose $N=N(T)$ so that 
\\
$$ N^{2(1-s)}\sim N^{-2+}N^{12(1-s)+}\Lambda(T)^{12+}$$
$$N(T) \sim \Lambda(T)^{\frac{12}{10s-8}+}.$$
This establishes Proposition 3 with
$$p(s)=\frac{2\cdot 12(1-s)}{10s-8}<2$$
\\
thus for $s>10/11$. We emphasize that for $1>s>10/11$ the second term in the conservation of the modified energy formula
 produces a smaller correction.
\end{proof}
{ \bf Acknowledgments}. I would like to thank Jim Colliander for introducing me to this problem. I am also grateful to him,
 for the beautiful discussions that we had while I was preparing the manuscript. 


\begin{thebibliography}{99}
\bibitem{bo} J. Bourgain, 
 {\it Refinements of Strichartz's inequality and applications to 2D NLS with critical nonlinearity,} Intern. Mat. Res. Notices,
 5, (1998), 253-283. 
\bibitem{bw} J. Bourgain, W. Wang, {\it Construction of blow-up solutions for the nonlinear Schr\"odinger equation with critical
 nonlinearity,} Ann. Scuola Norm. Sup. Pisa Cl. Sci. (4) 25 (1-2), (1997), 197-215. 
\bibitem{tc} T. Cazenave, {\it Semilinear 
Schr\"odinger equation,} AMS and C. I. M. S. Lecture Notes, 10, 2003. 
\bibitem{cw} T. Cazenave, F. B. Weissler {\it The Cauchy problem for the critical nonlinear 
Schr\"odinger equation in $H^{s}$,} Nonlinear Analysis Theory and Applications. 14, no.10 (1990), 807-836. 
\bibitem{ckstt1} J. Colliander, M. Keel, G. Staffilani, H. Takaoka, T. Tao, 
 {\it Global Well-posedness for Schr\"odinger equations with derivative,} SIAM J. Math. Anal. 33, no. 3, (2001), 649-669. 
\bibitem{ckstt2} J. Colliander, M. Keel, G. Staffilani, H. Takaoka, T. Tao, 
 {\it Sharp multi-linear periodic KdV estimates and applications,} J. Func. Anal. 211 (2004), 173-218.
\bibitem{ckstt3} J. Colliander, M. Keel, G. Staffilani, H. Takaoka, T. Tao, 
 {\it A refined Global Well-posedness result for Schr\"odinger equations with derivative,} SIAM J. Math. Anal. 34 (2002), 64-86.
\bibitem{ckstt4} J. Colliander, M. Keel, G. Staffilani, H. Takaoka, T. Tao, 
 {\it Sharp global well-posedness for KdV and modified KdV on $\Bbb R$  and
$\Bbb T$}, J. Amer. Math. Soc. 16 (2003), 705-749.
\bibitem{ckstt5} J. Colliander, M. Keel, G. Staffilani, H. Takaoka, T. Tao, 
 {\it Almost conservation laws and global rough solutions to a nonlinear Schr\"odinger
equation}, Math. Research Letters 9 (2002), 659-682.
\bibitem{ckstt6} J. Colliander, M. Keel, G. Staffilani, H. Takaoka, T. Tao, 
 {\it Global well-posedness and scattering for the energy-critical nonlinear Schr\"odinger
equation in ${\bf R^{3}}$}, to appear in the Annals of Mathematics.
\bibitem{crsw} J. Colliander, S. Raynor, C. Sulem, and J. D. Wright 
 {\it Ground state mass concentration in the $L^{2}$ critical nonlinear Schr\"odinger equation below $H^{1}$}, preprint.
\bibitem{gtv} J. Ginibre, Y. Tsutsumi and G. Velo, { \it On the Cauchy
problem for the Zakharov system,} J. Func. Anal. 151 (1997), 384-436.
\bibitem{kpv} C. E. Kenig, G. Ponce, L. Vega, {\it A bilinear estimate with applications to the Kdv equations,} 
J. Amer. Math. Soc 9 (1996), 573-603.
\bibitem{lpss} M. J. Landman, G. C. Papanicolaou, C. Sulem, and P-L. Sulem, {\it Rate of blow up for solutions
 of the nonlinear Schr\"odinger equation at critical dimensions,} Phys. Rev A. (3) 38 no. 8, (1988), 3837-3843.
\bibitem{l1} P.-L Lions, {\it The concentration compactness principle in the calculus of variations. The locally compact case. I,}
 Ann. Inst. Henri Poincar\'{e} Anal. Non Lin\'{e}aire 1, (1984), 109-145.
\bibitem{l2}P.-L Lions, {\it The concentration compactness principle in the calculus of variations. The locally compact case. II,}
 Ann. Inst. Henri Poincar\'{e} Anal. Non Lin\'{e}aire 1, (1984), 223-283.
\bibitem{mm} Y. Martel, F. Merle, {\it Blow up in finite time and dynamics of blow up solutions for the $L^{2}$-critical
 generalized KdV equation,} J. Amer. Math. Soc. 15 (2002), 617-664.
\bibitem{mr} F. Merle, P. Raphael, {\it Sharp upper bound on the blow-up rate 
for the critical nonlinear Schr\"odinger equation,} Geom. Funct. Anal. 13 (2003), 591-642.
\bibitem{mr1} F. Merle, P. Raphael, {\it Blow-up dynamic and upper bound on the blow-up rate 
for critical nonlinear Schr\"odinger equation,} to appear in Ann. Math.
\bibitem{mt} F. Merle, Y. Tsutsumi, {\it $L^{2}$ concentration of blow-up solutions for nonlinear Schr\"odinger equation
 with critical power nonlinearity,} J. Diff. Equations 84 (1990), 205-214.
\bibitem{mv} F. Merle, L. Vega, {\it Compactness at blow-up timw for $L^{2}$ solutions of the critical nonlinear Schr\"odinger equation
 in 2d,} Intern. Mat. Res. Notices, 8, (1998), 399-425.
\bibitem{hn} H. Nawa, {\it ``Mass concentration'' phenomenon for the nonlinear schr\"odinger equation with critical
 power nonlinearity,} Funkcial. Ekvac. 35 (1992), no. 1, 1-18.
\bibitem{hp} H. Pecher, {\it Global solutions with infinite energy for the 1-dimensional Zakharov system,} preprint
\bibitem{pg} G. Perelman {\it On the formation of singularities in solutions of the critical nonlinear Schr\"odinger equation,} 
Ann. Henri Poincar\'{e} 2 (4) (2001), 605-673.
\bibitem{ws} W. Strauss, {\it Existence of solitary waves in higher dimensions,} Commun. Math. Phys., 55 (1977), 149-162.
\bibitem{ss} C. Sulem, P-L. Sulem {\it The nonlinear Schr\"odinger equation: Self-focusing and wave collapse,} Series in Applied Mathematical Sciences, 
Springer, 139 (1999). 
\bibitem{tz} N. Tzirakis, {\it The Cauchy problem for the semilinear quintic Schr\"odinger equation in one dimension,} Differential Integral Equations,
 18 (2005), no 8, 947-960.
\bibitem{mw} M. Weinstein, {\it On the structure and formation of singularities of solutions to nonlinear dispersive evolution
 equations,} Comm. Partial Differential Equations, 11 (1986), 545-565.
\bibitem{mw1} M. Weinstein, {\it Nonlinear Schr\"odinger equations and sharp interpolation estimates,} Commun. Math. Phys., 87 (1983), 567-576.
\end{thebibliography}
\end{document}